\documentclass{article}
\usepackage{bm}
\usepackage{graphicx}
\usepackage{amsfonts}
\usepackage{romp}

\title{
Simple fractal calculus from fractal arithmetic}
\author{Diederik Aerts\\
Centrum Leo Apostel, Vrije Universiteit Brussel\\ 
Krijgskundestraat 33, 1160 Brussels, Belgium\\
diraerts@vub.ac.be\\[2ex]
Marek Czachor\\
Katedra Fizyki Teoretycznej i Informatyki Kwantowej, Politechnika Gda\'nska \\
ul. G. Narutowicza 11/12, 80-233 Gda\'nsk, Poland\\
mczachor@pg.edu.pl\\[2ex]
Maciej Kuna\\
Katedra Rachunku Prawdopodobie\'nstwa i Biomatematyki, Politechnika Gda\'nska \\
ul. G. Narutowicza 11/12, 80-233 Gda\'nsk, Poland\\
maciek@mif.pg.gda.pl}
\begin{document}
\maketitle

\begin{abstract}
Non-Newtonian calculus that starts with elementary non-Diophantine arithmetic operations of a Burgin type is applicable to all fractals whose cardinality is continuum. The resulting definitions of derivatives and integrals are simpler from what one finds in the more traditional literature of the subject, and they often work in the cases where the standard methods fail. As an illustration, we perform a Fourier transform of a real-valued function with Sierpi\'nski-set domain. The resulting formalism is as simple as the usual undergraduate calculus.
\end{abstract}


\section{Introduction}

Apparently, the first attempt of a Fourier-type analysis on fractals can be found in studies of diffusion on Sierpi\'nski gaskets \cite{Kusuoka,Goldstein}. A generator of the diffusion process plays there the same role as a Laplacian on a manifold, so the corresponding eigenfunction expansion may be regarded as a form of harmonic analysis.  An alternative route to eigenfunction expansions on fractals is to define Laplacians or gradients more directly. Here certain approaches begin with Dirichlet forms on self-similar fractals, or one takes as a departure point discrete Laplacians and performes an appropriate limit \cite{Kusuoka1989,Kigami1989,Kigami,Strichartz}. Four alternative definitions of a gradient (due to Kusuoka, Kigami, Strichartz and Teplyaev) are discussed in this context in \cite{Tep}. Self-similarity is typically an important technical assumption. Although Laplacians defined in the above ways cannot be regarded as second-order operators, an approach where Laplacians are indeed second-order is nevertheless possible and was introduced by Fujita \cite{Fujita1987,Fujita} and further developed by a number of authors \cite{FZ,F,Z,Arzt,Kess}.

A second traditional approach to harmonic analysis on fractals comes from the notion of self-similar fractal measures. The classic result of Jorgensen and Pedersen \cite{JP} states that the method works for certain fractals, such as the quaternary Cantor set, but fails in the important case of the ternary middle-third Cantor set. Moreover, the method is inapplicable in realistic cases of non-self-similar fractals.

Quite recently we have shown \cite{ACK2} that a `non-Newtonian' calculus, based on Burgin's non-Diophantine arithmetic, \cite{MC,Burgin1,Burgin2,Burgin3,ACK1} leads to a simple and very efficient construction of a Fourier transform on fractals of a Cantor type. Gradients and Laplacians are here, respectively, first- and second-order differential operators, and self-similarity plays no role whatsoever. There is completely no difficulty with Fourier analysis of functions mapping arbitrary Cantor sets into themselves, so Jorgensen-Pedersen-type restrictions are no longer valid. 
The question of Fourier analysis on fractals is important for the problem of momentum representation in quantum mechanics on fractal space-times. 
Another recent application of the calculus is deformation quantization with minimal length  \cite{Blaszak}, and the problem of wave equations on space-times modeled by Cartesian products of different fractals  \cite{MC2017}.

The goal of the present paper is to show explicitly how to apply the non-Diophantine framework to fractals more general than the Cantor set. We explicitly perform the construction for a double cover of a Sierpi\'nski set. Similarly to the Cantor case, self-similarity is inessential. What is important, however, is the existence of a bijection $f$ between the fractal in question and $\mathbb{R}$.

In the Sierpi\'nski case the bijection has a space-filling property reminiscent of Peano curves \cite{Peano}. The very idea that there are links between Sierpi\'nski-type fractals and space-filling curves is not new, and was used by Molitor {\it et al.\/} in \cite{Peano1} in their construction of Laplacians on fractals. However, in all other respects the approach from \cite{Peano1} is different from what we discuss below. The idea of employing a one-dimensional integration for finding higher dimensional integrals is known \cite{Wiener}, but apparently has not been used in fractal contexts so far.

Since any fractal whose cardinality is continuum can be equipped with a bijection mapping it into $\mathbb{R}$, the construction is quite universal. From a practical perspective, the only difficulty is to find the bijection explicitly, but once we have found it the remaining procedure is systematic and easy to work with.

Here, out of a multitude of possible illustrations of the formalism we have decided to discuss the case of a sine Fourier transform of a real-valued function with Sierpi\'nski-set domain. One can directly judge applicability of the method by visually inspecting the quality of the resulting finite-term reconstruction of the signal.

\section{Sierpi\'nski set}

Consider $x\in\mathbb{R}_+$ and its ternary representation $x=(t_n\dots t_0.t_{-1}t_{-2}\dots)_3$. If $x$ has two different ternary representations, we choose the one that ends with infinitely many 2s. All finite-digit numbers are thus represented by infinite sequences, which is perhaps unusual but reduces ambiguity of the inverse algorithm, as we shall see shortly.  Keeping the digits unchanged let us change the base from 3 to 4, i.e.
\begin{eqnarray}
x=(t_n\dots t_0.t_{-1}t_{-2}\dots)_3\mapsto (t_n\dots t_0.t_{-1}t_{-2}\dots)_4=y.
\end{eqnarray}
The quaternary representation of $y$ is unique, and it does not involve the digit 3. Next, let us parametrize the quaternary digits in a binary way, but written in a column form: $0=\begin{array}{c}0\\0\end{array}$, $1=\begin{array}{c}0\\1\end{array}$, $2=\begin{array}{c}1\\0\end{array}$,
$3=\begin{array}{c}1\\1\end{array}$. So, $y$ has been converted into a pair of binary sequences,
\begin{eqnarray}
(t_n\dots t_0.t_{-1}t_{-2}\dots)_4
\mapsto
\left(
\begin{array}{c}
a_n\dots a_0.a_{-1}a_{-2}\dots\\
b_n\dots b_0.b_{-1}b_{-2}\dots
\end{array}
\right)_2
\end{eqnarray}
where $(a_j,b_j)\neq (1,1)$ for any $j$. The resulting sequences are in a one-one relation with the $x$ we have started with. Each of the two sequences defines a number in binary notation: we have mapped $x$ into a point of the plane $\mathbb{R}_+\times \mathbb{R}_+$. The image of $\mathbb{R}_+$ under our algorithm defines a Sierpi\'nski-type set. The algorithm is not invertible. Indeed, take the point $(1,1)$. Depending on the way we represent it binarily we find
\begin{eqnarray}
\left(
\begin{array}{c}
1.(0)\\
0.(1)
\end{array}
\right)_2
&\mapsto&
\big(2.(1)\big)_4
\mapsto
\big(2.(1)\big)_3=2.5,
\nonumber
\end{eqnarray}
and
\begin{eqnarray}
\left(
\begin{array}{c}
0.(1)\\
1.(0)
\end{array}
\right)_2
&\mapsto&
\big(1.(2)\big)_4
\mapsto
\big(1.(2)\big)_3=2.\nonumber
\end{eqnarray}
The ambiguity comes from the two identifications
\begin{eqnarray}
(1,1)=
\left(
\begin{array}{c}
1.(0)\\
0.(1)
\end{array}
\right)_2,
\quad
(1,1)=
\left(
\begin{array}{c}
0.(1)\\
1.(0)
\end{array}
\right)_2.
\end{eqnarray}
However, if we write the above two relations as
\begin{eqnarray}
(1,1)_-
&=&
\left(
\begin{array}{c}
1.(0)\\
0.(1)
\end{array}
\right)_2,
\\
(1,1)_+
&=&
\left(
\begin{array}{c}
0.(1)\\
1.(0)
\end{array}
\right)_2,
\end{eqnarray}
and treat the two points $(1,1)_\pm$ as belonging to two different sides of an oriented plane, the ambiguity of the inverse alorithm disappears. The relation
\begin{eqnarray}
(1,1)_-
\leftrightarrow
2.5,
\quad
(1,1)_+
\leftrightarrow
2
\end{eqnarray}
is one-one.

Let us therefore index with `$+$' (resp. `$-$') those pairs $(a,b)$ where $a$ is a rational number represented by a binary sequence involving $(1)_2$ (resp. $(0)_2$), and $b$ is a rational number whose binary representation contains $(0)_2$ (resp. $(1)_2$).
In both cases $(a_j,b_j)\neq (1,1)$ by construction. The corresponding rational numbers involve, respectively, $(2)_3$ and $(1)_3$.

But what about the other cases, such as $a$, $b$ irrational, or $a$ irrational but $b$ rational? It turns out that the ambiguity is absent.
In order to prove it, first of all note that we did not have to consider the cases
\begin{eqnarray}
(1,1)=
\left(
\begin{array}{c}
0.(1)\\
0.(1)
\end{array}
\right)_2,
\quad
(1,1)=
\left(
\begin{array}{c}
1.(0)\\
1.(0)
\end{array}
\right)_2,\label{(00)}
\end{eqnarray}
since the pairs $(a_j,b_j)=(1,1)$ cannot appear as a result of the algorithm, and two infinite sequences of 0s would imply that
$(t_n\dots t_0.t_{-1}t_{-2}\dots)_4$ ends with an infinite sequence of 0s, a form excluded by the algorithm.

The same mechanism eliminates all the remaining ambiguities:

(A) If $a$, $b$ are both irrational, or $a$ is irrational and $b$ rational-periodic, their binary forms are unique.

(B) If $a$ is irrational (or rational-periodic), but $b$ rational non-periodic, then $b$ cannot end with infinitely many 1s, as it would mean that $a$ ends with infinitely many 0s. So these cases are again unique. Conclusions of (A) and (B) are unchanged if one interchanges the roles of $a$ and $b$.

(C) The only ambiguity appears if $a$ ends with infinitely many 0s, but $b$ with infinitely many 1s (or the other way around). But this is the case we have started with.

In cases (A) and (B), we identify $(a,b)_+=(a,b)_-=(a,b)$. Only the (countable) case (C) requires a two-sided plane $(a,b)_+\neq (a,b)_-$. The case (C) occurs for those $x\in\mathbb{R}$ whose ternary representation ends with $(2)_3$ or $(1)_3$. Only the latter numbers are mapped into $(a,b)_-$.

As we can see, what we have constructed is a version of a double cover of the Sierpi\'nski set.

Our algorithm defines an injective map $g_+$ of $\mathbb{R}_+$ into a two-sided plane, with the above-mentioned identifications. Let us extend $g_+$ to $g$ by $g(|x|)=g_+(|x|)$, $g(-|x|)=-g_+(|x|)$. The image $S=g(\mathbb{R})$ is our definition of the Sierpi\'nski set. Denoting $f=g^{-1}$, $f:S\to\mathbb{R}$ we obtain Burgin's arithmetic intrinsic to $S$,
\begin{eqnarray}
x\oplus y &=& f^{-1}\big(f(x)+f(y)\big),\\
x\ominus y &=& f^{-1}\big(f(x)-f(y)\big),\\
x\odot y &=& f^{-1}\big(f(x)f(y)\big),\\
x\oslash y &=& f^{-1}\big(f(x)/f(y)\big).
\end{eqnarray}
Fig.~1 shows the set $f^{-1}\big([0,1)\big)$. The set is self-similar and its Hausdorff dimension is $\log_2 3$.
\begin{figure}
\includegraphics[width=8 cm]{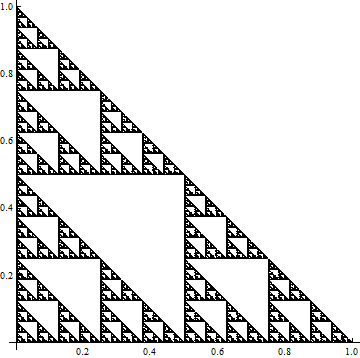}
\caption{The inverse image $f^{-1}([0,1))$, view from the positive side of the oriented plane.}
\end{figure}

\section{Arithmetic on the Sierpi\'nski set}

Let us begin with neutral elements of addition and multiplication in $S$. By definition, $0'\oplus x=x$, $1'\odot x=x$, where
\begin{eqnarray}
0' &=&f^{-1}(0)=(0,0)\in S,\\
1' &=&f^{-1}(1)=(1,0)_+\in S.
\end{eqnarray}

In a non-Diophantine arithmetic, multiplication is a repeated addition in the following sense \cite{ACK2}. Let $n\in\mathbb{N}$ and $n'=f^{-1}(n)\in S$. Then
\begin{eqnarray}
n'\oplus m'
&=&
(n+m)',\\
n'\odot m'
&=&
(nm)'\\
&=&
\underbrace{m'\oplus\dots \oplus m'}_{n\rm{ times}}.
\end{eqnarray}

A power function $A(x)=x\odot\dots \odot x$ ($n$ times) is denoted by $x^{n'}$, which is consistent with
\begin{eqnarray}
x^{n'}\odot x^{m'}=x^{(n+m)'}=x^{n'\oplus m'}.
\end{eqnarray}

All integers are here represented by pairs of integers, a representation somewhat similar to complex numbers, but with different rules of addition and multiplication, as illustrated by
\begin{eqnarray}
3'\oplus 4' = (2,0)_+\oplus (1,2)_+=7'=(3,0)_+
\end{eqnarray}
(numbers represented in decimal form).
\begin{figure}
\includegraphics[width=8 cm]{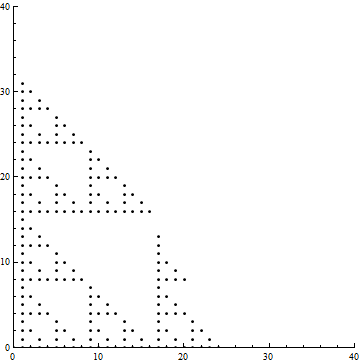}
\caption{The image of the first 200 natural numbers, $f^{-1}(\{1,\dots, 200\})$. All natural numbers are mapped into the positive side of the oriented plane.}
\end{figure}

\section{Calculus}

The map $g$ is not continuous (in Euclidean metric topology of $\mathbb{R}^2$), as illustrated by the following generic example. Consider
\begin{eqnarray}
g(1+1/3^n)
&=&
g\big(1.\underbrace{0\dots 0}_n(2)_3\big)=\left(
\begin{array}{c}
0.0\dots 0(1)\\
1.0\dots 0(0)
\end{array}
\right)_2,\nonumber\\
g(1-1/3^n)
&=&
g\big(0.\underbrace{2\dots 21}_n(2)_3\big)=\left(
\begin{array}{c}
0.1\dots 10(1)\\
0.0\dots 01(0)
\end{array}
\right)_2.\nonumber
\end{eqnarray}
The function is discontinuous at $x=1$,
\begin{eqnarray}
\lim_{x\to 1_+}g(x)=(0,1)_+,\quad \lim_{x\to 1_-}g(x) =(1,0)_+,
\nonumber
\end{eqnarray}
but the argument will not work at $x=0$,
\begin{eqnarray}
\lim_{x\to 0_+}g(x)=\lim_{x\to 0_-}g(x)=(0,0)=0',
\end{eqnarray}
so the limit $h\to 0'$ is unambiguous.
Let us recall that $0'$ is the neutral element of Burgin's non-Diophantine addition.

The derivative of a function $A:S\to S$ is defined by \cite{MC,ACK1,ACK2}
\begin{eqnarray}
\frac{DA(x)}{Dx}
&= &
\lim_{h\to 0'}\Big(A(x\oplus h)\ominus A(x)\Big)\oslash h.
\end{eqnarray}
A Laplacian on $S$ is just the second derivative
\begin{eqnarray}
\Delta A(x)=\frac{D}{Dx}\frac{DA(x)}{Dx}.
\end{eqnarray}
Our $\Delta$ differs from the other definitions of Laplacians on Sierpi\'nski sets occurring in the literature \cite{Kigami,Strichartz}, but is simple and easy to work with. The formalism from \cite{ACK2,MC,ACK1} can be applied here with no modification, including integration, complex numbers, harmonic analysis, differential equations and so on.

What was {\it not\/} explained in \cite{ACK2,MC,ACK1} was how to proceed with functions that do not map the fractal in question into itself.
So, consider two sets, $\mathbb{X}$ and $\mathbb{Y}$ say, equipped with bijections $f_\mathbb{Y}:\mathbb{Y}\to \mathbb{R}$ and $f_{\mathbb{X}}:\mathbb{X}\to \mathbb{R}$, and arithmetics $\{\oplus_{\mathbb{Y}},\odot_{\mathbb{Y}}:\mathbb{Y}\times \mathbb{Y}\to \mathbb{Y}\}$,  $\{\oplus_{\mathbb{X}},\odot_{\mathbb{X}}:\mathbb{X}\times \mathbb{X}\to \mathbb{X}\}$, defined by $f_{\mathbb{Y}}$ and $f_{\mathbb{X}}$.
The bijection $f=f_{\mathbb{Y}}^{-1}\circ f_{\mathbb{X}}: \mathbb{X}\to \mathbb{Y}$ makes it possible to consider derivatives of functions $A:\mathbb{X}\to \mathbb{Y}$. Let $0'_{\mathbb{X}}=f^{-1}_{\mathbb{X}}(0)$ be the neutral element of addition in $\mathbb{X}$, and $f(0'_{\mathbb{X}})=0'_{\mathbb{Y}}$ the one in $\mathbb{Y}$.
We define
\begin{eqnarray}
\frac{DA(x)}{Dx}
&=&
\lim_{h\to 0'_{\mathbb{X}}}\Big(A(x\oplus_{\mathbb{X}} h)\ominus_{\mathbb{Y}} A(x)\Big)\oslash_{\mathbb{Y}} f(h)\label{deriv}\\
&=&
\lim_{h\to 0}\Big(A\big(x\oplus_{\mathbb{X}} f^{-1}_\mathbb{X}(h)\big)\ominus_{\mathbb{Y}} A(x)\Big)
\oslash_{\mathbb{Y}} f^{-1}_\mathbb{Y}(h)\label{deriv0}\\
&=&
f_{\mathbb{Y}}^{-1}\left(\frac{d}{df_{\mathbb{X}}(x)}\underbrace{f_{\mathbb{Y}}\circ A\circ f_{\mathbb{X}}^{-1}}_a\big[f_{\mathbb{X}}(x)\big]\right).\label{deriv1}
\end{eqnarray}
Details of the transition from (\ref{deriv})--(\ref{deriv0}) to (\ref{deriv1}) can be found in the Appendix.
Accordingly,
\begin{eqnarray}
\frac{DA(x)}{Dx}=f_{\mathbb{Y}}^{-1}\circ a'\circ f_{\mathbb{X}}(x)
\end{eqnarray}
where $a'(x)=\lim_{h\to 0}\big(a(x+h)-a(x)\big)/h$.

The integral is defined in a way guaranteeing the fundamental laws of calculus, relating derivatives and integrals.
So, let $A:\mathbb{X}\to \mathbb{Y}$, and
\begin{eqnarray}
\int_Y^X A(x)Dx
&=&
f^{-1}_{\mathbb{Y}}\left(\int_{f_{\mathbb{X}}(Y)}^{f_{\mathbb{X}}(X)}f_{\mathbb{Y}}\circ A\circ f_{\mathbb{X}}^{-1}(x)dx\right)
\end{eqnarray}
where $\int a(x)dx$ is the usual (say, Lebesgue) integral of a function $a:\mathbb{R}\to \mathbb{R}$.

In the Appendix we prove that
\begin{eqnarray}
\frac{D}{DX}\int_Y^X A(x)Dx
&=&A(X),\label{int1}\\
\int_Y^X\frac{DA(x)}{Dx}Dx
&=&
A(X)\ominus_{\mathbb{Y}} A(Y)\label{int2}.
\end{eqnarray}
In the next section we apply the formalism to the problem of harmonic analysis on Sierpi\'nski sets.

\section{Example: Fourier transform on $S$}

Let $\mathbb{X}=S$ and $\mathbb{Y}=\mathbb{R}$. Consider the function $A:S\to\mathbb{R}$ (Fig.~3),
\begin{eqnarray}
A(x)
&=&
\left\{
\begin{array}{cl}
1 &\textrm{for }x\in f_S^{-1}\big((0,1)\big)\\
-1 &\textrm{for }x\in f_S^{-1}\big((-1,0)\big)\\
0 & \textrm{otherwise}
\end{array}
\right.\label{A(x)}
\end{eqnarray}
Since $S\cap(-1,0)^2=f_S^{-1}\big((-1,0)\big)$,  $S\cap(0,1)^2=f_S^{-1}\big((0,1)\big)$, $S\cap \{(0,0)\}=f_S^{-1}\big(0\big)$,
we introduce $a:\mathbb{R}\to \mathbb{R}$,
\begin{eqnarray}
a(x)
&=&
\left\{
\begin{array}{cl}
1 &\textrm{for }x\in (0,1)\\
-1 &\textrm{for }x\in (-1,0)\\
0 & \textrm{otherwise}
\end{array}
\right.
\end{eqnarray}
Employing $f_{\mathbb{X}}=f_S$, $f_{\mathbb{Y}}=\textrm{id}_{\mathbb{R}}$, we get
\begin{eqnarray}
a=f_{\mathbb{Y}}\circ A\circ f_{\mathbb{X}}^{-1}= A\circ f_{S}^{-1}.
\end{eqnarray}
In order to perform Fourier analysis of $A$ we have to introduce the basis of sines and cosines along the lines of \cite{ACK2}, but adapted to the present context. The scalar product of two functions $G_j:S\to\mathbb{R}$, $G_j=g_j\circ f_S$, $g_j:\mathbb{R}\to \mathbb{R}$, $j=1,2$, reads
\begin{eqnarray}
\langle G_1|G_2\rangle
&=&
\int_{\ominus_{\mathbb{Y}} T}^T G_1(x)\odot_{\mathbb{Y}}G_2(x)Dx\\
&=&
f_{\mathbb{Y}}^{-1}\left(\int_{-f_S(T)}^{f_S(T)} g_1(x)g_2(x)dx\right)
\\
&=&
\int_{-1}^{1} g_1(x)g_2(x)dx=\langle g_1|g_2\rangle\label{scal}
\end{eqnarray}
where $\ominus_S T=0'_S\ominus_ST=f^{-1}_S(-1)$. In our case $T=1'_S=f^{-1}_S(1)$ (the neutral element of multiplication in $S$).
The fact that $\langle g_1|g_2\rangle$ and $\langle G_1|G_2\rangle$ involve the same symbol of the scalar product will not lead to ambiguities.

Denoting
\begin{eqnarray}
c_n(y)
&=&
\cos n\pi y,\quad n>0\\
s_n(y)
&=&
\sin n\pi y,\quad n>0\\
c_0(y)
&=&
1/\sqrt{2},\\
s_0(y)
&=&
0,\\
C_n(x) &=& c_n\big(f_S(x)\big),\\
S_n(x) &=& s_n\big(f_S(x)\big),
\end{eqnarray}
we can apply the standard resolution of unity,
\begin{eqnarray}
\delta(x-y)
&=&
\sum_{n\geq 0}\Big(c_n(x)c_n(y)+s_n(x)s_n(y)\Big),
\end{eqnarray}
and finally obtain
\begin{eqnarray}
A(x)
&=&
\sum_{n\geq 0}
\Big(
C_n(x)\langle C_n|A\rangle
+
S_n(x)\langle S_n|A\rangle
\Big)
\\
&=&
\sum_{n\geq 0}
\Big(
C_n(x)\langle c_n|a\rangle
+
S_n(x)\langle s_n|a\rangle
\Big)\\
&=&
\sum_{n> 0}
\frac{2 \left(1-(-1)^n\right)}{n \pi }S_n(x)\label{An}
\end{eqnarray}
\begin{figure}
\includegraphics[width=8 cm]{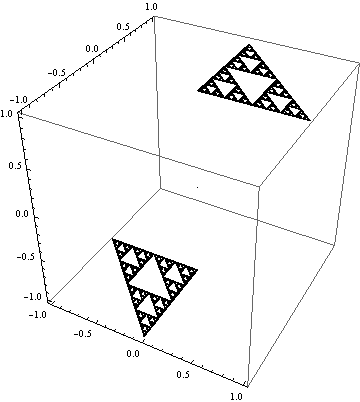}
\caption{Function $A(x)$ defined in (\ref{A(x)}).}
\end{figure}
\begin{figure}
\includegraphics[width=8 cm]{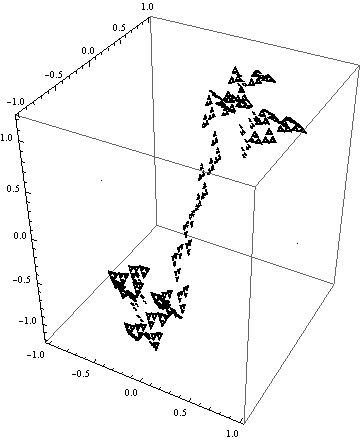}
\caption{Finite-sum reconstruction of $A(x)$ with 5 Fourier terms in (\ref{An}) ...}
\end{figure}
\begin{figure}
\includegraphics[width=8 cm]{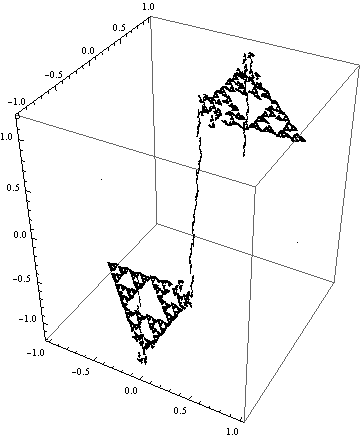}
\caption{...and with 50 terms. The Gibbs phenomenon is clearly visible.}
\end{figure}
Figures 4 and 5 illustrate finite-sum Fourier reconstructions of the function plotted in Fig.~3. The negative-side of $S$ occurs only as the image of those $x\in\mathbb{R}$ whose ternary representation ends with infinitely many 1s, whereas all finite-ternary-digit numbers are mapped into the positive side of $S$. This leads to the practical moral: Plotting functions with domains in $S$ we can concentrate exclusively on the positive side of $S$, unless one employs a symbolic algorithm that recognizes numbers ending with $(1)_3$.

\section{Conclusion}

The formalism that starts with non-Diophantine addition and multiplication is applicable to a general class of fractals, including those of a Sierpi\'nski type. For a fractal whose cardinality equals continuum the existence of a  bijection $f$ follows directly from the fact that the cardinality of $\mathbb{R}$ is the same, so the resulting paradigm is universal. From the point of view of applications the only nontrivial element of the construction is to find an explicit $f$. Once we have found it, the remaining procedures are systematic and as simple as an undergraduate calculus.

\section*{Appendix}

Here we give detailed proofs of formulas (\ref{deriv}), (\ref{int1}), and (\ref{int2}).

One begins with the limit $h\to 0'_{\mathbb{X}}$ which is defined as follows. Consider the commutative diagram
\begin{eqnarray}
\begin{array}{rcl}
\mathbb{X}                & \stackrel{A}{\longrightarrow}       & \mathbb{Y}               \\
f_\mathbb{X}{\Big\downarrow}   &                                     & {\Big\downarrow} f_\mathbb{Y}  \\
\mathbb{R}                & \stackrel{a}{\longrightarrow}   & \mathbb{R}  
\end{array}
\end{eqnarray}
In our example we had $\mathbb{Y}=\mathbb{R}$ and $f_\mathbb{Y}=\textrm{id}_\mathbb{R}$, so $f^{-1}_\mathbb{Y}$ is trivially continuous, which justifies the transition from (\ref{51}) to (\ref{52}) below. However, the formalism works also in cases where $f^{-1}_\mathbb{Y}$  is discontinuous, but the limit is understood as 
\begin{eqnarray}
\lim_{X\to X_0}A(X) 
&=&
f_\mathbb{Y}^{-1}\Big(\lim_{x\to f_\mathbb{X}(X_0)}a(x)\Big),
\end{eqnarray}
for $X_0\in\mathbb{X}$. 
The latter is logically equivalent to {\it defining\/} the derivative directly by the end result of the calculation we give below, that is
\begin{eqnarray}
\frac{DA(x)}{Dx}
&:=&
f_{\mathbb{Y}}^{-1}\left(\frac{da\big[f_{\mathbb{X}}(x)\big]}{df_{\mathbb{X}}(x)}\right).
\end{eqnarray}

The following sequence of transformations is instructive, and it explains why differentiability of the bijections is unnecessary to make the calculation work:
\begin{eqnarray}
\frac{DA(x)}{Dx}
&=&
\lim_{h\to 0'_{\mathbb{X}}}\Big(A(x\oplus_{\mathbb{X}} h)\ominus_{\mathbb{Y}} A(x)\Big)\oslash_{\mathbb{Y}} f(h)\\
&=&
\lim_{h\to 0'_{\mathbb{X}}}f_{\mathbb{Y}}^{-1}\left(\frac{f_{\mathbb{Y}}\big(A(x\oplus_{\mathbb{X}} h)\ominus_{\mathbb{Y}} A(x)\big)}{f_{\mathbb{Y}}\big( f(h)\big)}\right)\\
&=&
\lim_{h\to 0'_{\mathbb{X}}}f_{\mathbb{Y}}^{-1}\left(\frac{f_{\mathbb{Y}}\big(A(x\oplus_{\mathbb{X}} h)\ominus_{\mathbb{Y}} A(x)\big)}{f_{\mathbb{X}}(h)}\right)\\
&=&
\lim_{h\to 0'_{\mathbb{X}}}f_{\mathbb{Y}}^{-1}\left(\frac{f_{\mathbb{Y}}\big(A(x\oplus_{\mathbb{X}} h)\big)-f_{\mathbb{Y}}\big( A(x)\big)}{f_{\mathbb{X}}(h)}\right)\\
&=&
\lim_{h\to 0'_{\mathbb{X}}}f_{\mathbb{Y}}^{-1}\left(\frac{f_{\mathbb{Y}}\Bigg(A\Big(f_{\mathbb{X}}^{-1}\big[f_{\mathbb{X}}(x)+ f_{\mathbb{X}}(h)\big]\Big)\Bigg)-f_{\mathbb{Y}}\Bigg(A\Big(f_{\mathbb{X}}^{-1}\big[f_{\mathbb{X}}(x)\big]\Big)\Bigg)}{f_{\mathbb{X}}(h)}\right)\\
&=&
\lim_{h\to 0'_{\mathbb{X}}}f_{\mathbb{Y}}^{-1}\left(\frac{f_{\mathbb{Y}}\circ A\circ f_{\mathbb{X}}^{-1}\big[f_{\mathbb{X}}(x)+f_{\mathbb{X}}(h)\big]-f_{\mathbb{Y}}\circ A\circ f_{\mathbb{X}}^{-1}\big[f_{\mathbb{X}}(x)\big]}{f_{\mathbb{X}}(h)}\right)\label{51}\\
&=&
f_{\mathbb{Y}}^{-1}\left(\lim_{h\to 0}\frac{f_{\mathbb{Y}}\circ A\circ f_{\mathbb{X}}^{-1}\big[f_{\mathbb{X}}(x)+h\big]-f_{\mathbb{Y}}\circ A\circ f_{\mathbb{X}}^{-1}\big[f_{\mathbb{X}}(x)\big]}{h}\right)\label{52}\\
&=&
f_{\mathbb{Y}}^{-1}\left(\frac{d}{df_{\mathbb{X}}(x)}\underbrace{f_{\mathbb{Y}}\circ A\circ f_{\mathbb{X}}^{-1}}_a\big[f_{\mathbb{X}}(x)\big]\right)\\
&=&
f_{\mathbb{Y}}^{-1}\left(\frac{da\big[f_{\mathbb{X}}(x)\big]}{df_{\mathbb{X}}(x)}\right).
\end{eqnarray}
This proves (\ref{deriv1}). 
Now, let $A:\mathbb{X}\to \mathbb{Y}$, and define
\begin{eqnarray}
\int_Y^X A(x)Dx
&=&
f^{-1}_{\mathbb{Y}}\left(\int_{f_{\mathbb{X}}(Y)}^{f_{\mathbb{X}}(X)}a(x)dx\right)\\
&=&
f^{-1}_{\mathbb{Y}}\left(\int_{f_{\mathbb{X}}(Y)}^{f_{\mathbb{X}}(X)}f_{\mathbb{Y}}\circ A\circ f_{\mathbb{X}}^{-1}(x)dx\right).
\end{eqnarray}
Denoting $b(y)=\int_{f_{\mathbb{X}}(Y)}^{y}a(x)dx$ we rewrite
\begin{eqnarray}
\int_Y^X A(x)Dx
&=&
f^{-1}_{\mathbb{Y}}\Big(b\big(f_{\mathbb{X}}(X)\big)\Big),
\end{eqnarray}
and thus (\ref{deriv1}) implies
\begin{eqnarray}
\frac{D}{DX}\int_Y^X A(x)Dx
&=&
f^{-1}_{\mathbb{Y}}\left(\frac{db\big(f_{\mathbb{X}}(X)\big)}{df_{\mathbb{X}}(X)}\right)
\\
&=&
f^{-1}_{\mathbb{Y}}\Big(a\big(f_{\mathbb{X}}(X)\big)\Big)
\\
&=&
f^{-1}_{\mathbb{Y}}\circ a\circ f_{\mathbb{X}}(X)
\\
&=&
f^{-1}_{\mathbb{Y}}\circ f_{\mathbb{Y}}\circ A\circ f_{\mathbb{X}}^{-1}\circ f_{\mathbb{X}}(X)=A(X).
\end{eqnarray}
And the other way around,
\begin{eqnarray}
\int_Y^X\frac{DA(x)}{Dx}Dx
&=&
\int_Y^X f_{\mathbb{Y}}^{-1}\circ a'\circ f_{\mathbb{X}}(x)Dx\\
&=&
f^{-1}_{\mathbb{Y}}\left(\int_{f_{\mathbb{X}}(Y)}^{f_{\mathbb{X}}(X)}f_{\mathbb{Y}}\circ f_{\mathbb{Y}}^{-1}\circ a'\circ f_{\mathbb{X}}\circ f_{\mathbb{X}}^{-1}(x)dx\right)
\\
&=&
f^{-1}_{\mathbb{Y}}\left(\int_{f_{\mathbb{X}}(Y)}^{f_{\mathbb{X}}(X)}a'(x)dx\right)
\\
&=&
f^{-1}_{\mathbb{Y}}\Big(a\big(f_{\mathbb{X}}(X)\big)-a\big(f_{\mathbb{X}}(Y)\big)\Big)
\\
&=&
f^{-1}_{\mathbb{Y}}\Big(f_{\mathbb{Y}}\circ A\circ f_{\mathbb{X}}^{-1}\big(f_{\mathbb{X}}(X)\big)-f_{\mathbb{Y}}\circ A\circ f_{\mathbb{X}}^{-1}\big(f_{\mathbb{X}}(Y)\big)\Big)
\\
&=&
f^{-1}_{\mathbb{Y}}\Big(f_{\mathbb{Y}}\big(A(X)\big)-f_{\mathbb{Y}}\big(A(Y)\big)\Big)
\\
&=&
A(X)\ominus_{\mathbb{Y}} A(Y).
\end{eqnarray}

\end{document}